\documentclass[12pt]{amsart}
\usepackage{amsmath}
\usepackage{amssymb}
\usepackage{bm}
\usepackage{graphicx}
\usepackage{verbatim}
\newtheorem{Theorem}{Theorem}

\newtheorem{Corollary}{Corollary}

\title[Breakdown of Semiclassical Correspondence]{Eisenstein Series and Breakdown of Semiclassical Correspondence}

\author{Shimon Brooks}
\address{Department of Mathematics, Bar-Ilan University, Ramat-Gan, 5290002 Israel}
\email{brookss@math.biu.ac.il}

\thanks{The author was partially supported by Israel Science Foundation grant 1119/13 and a Marie-Curie Career Integration Grant.}

\begin{document}
\maketitle

{\em Abstract: }  We consider certain Lagrangian states associated to unstable horocycles on the modular surface $PSL(2,\mathbb{Z})\backslash\mathbb{H}$, and show that for sufficiently large logarithmic times, expectation values for the wave propagated states diverge from the classical transport along geodesics.  This is due to the fact that these states ``escape to the cusp" very quickly, at logarithmic times, while the geodesic flow continues to equidistribute on the surface.  The proof relies crucially on the analysis of expectation values for Eisenstein series initiated by Luo-Sarnak and Jakobson, based on subconvexity estimates for relevant $L$-functions--- in other words, this is a very special case in which we can analyze long time propagation explicitly with tools from arithmetic.

\section{Introduction}
A central problem in asymptotic analysis is the description of high-frequency wave propagation to long times.  For our purposes, let $M=PSL(2,\mathbb{Z})\backslash\mathbb{H}$ be the modular surface, and $\phi_r$ an $L^2$-normalized eigenfunction of the Laplacian satisfying 
$$(\Delta+\frac14+r^2)\phi_r=0$$
or, more generally, an approximate eigenfunction or {\bf quasimode} satisfying
$$\left\|(\Delta+\frac14+r^2)\phi_r\right\|_{L^2} \leq \omega(r) \|\phi_r\|_{L^2}$$
for some function $\omega(r)=o(r)$ as $r\to\infty$.  Since $|\phi_r|^2$ gives the probability distribution for position of the state $\phi_r$, we can form the measures 
$$\mu_{\phi_r}(g) = \langle g\phi_r, \phi_r\rangle$$
on $M$, giving the expectation value for $g$ of the state $\phi_r$.  This can be extended to all observables $g\in C^\infty(M)$ with a microlocal lift, but for the purposes of this paper we will content ourselves with observables $g\in C^\infty(M)$ depending only on position.

We are interested in the behavior of these expectation values under two transformations:  the wave flow $U^t$ acting on eigenfunctions via $\phi_r\mapsto e^{itr}\phi_r$, and the classical dynamics $f \mapsto (f\circ a_t)$, where $a_t=\begin{pmatrix} e^{t/2} & \\ & e^{-t/2} \end{pmatrix}$ is the time-one geodesic flow on $M$. 
For short time scales, the Egorov Theorem (see eg. \cite[\S11.4]{EvansZworski}) shows that the two expectations
\begin{equation}\label{sc c}
\big\langle g U^t\phi_r, U^t\phi_r\big\rangle \sim \big\langle g (\phi_r\circ a_t), (\phi_r\circ a_t)\big\rangle
\end{equation}
are asymptotic as $r\to\infty$.  It is known \cite{BouRob} that this {\bf semiclassical correspondence} is valid up to the Ehrenfest time $|t|\leq T_E=c\log{r}$ for some $c>0$.  This is proved by writing out an explicit Taylor expansion for the wave propagation, whose main term is given by the classical geodesic-flow; this method breaks down at the Ehrenfest time when the error terms in the expansion become larger than the main term.

However, numerical evidence (see eg. \cite{TomsovicHeller, TomsovicHeller2}) suggests that the semiclassical approximation could be valid for much longer times, perhaps even to powers of $r$.  Of particular interest in this regard is the evolution of Lagrangian states associated to unstable Lagrangian manifolds \cite{SchubertExpectations}, since the classical dynamics contract the stable direction, and thus the image under the geodesic flow tends to concentrate on unstable leaves.

Though one may expect the semiclassical correspondence (\ref{sc c}) of expectation values to hold for larger times on a generic manifold, the main result of this paper is that it cannot hold in full generality beyond logarithmic times.  Namely, on the modular surface we will construct Lagrangian states associated to unstable horocycles relative to the cusp, whose propagation under the wave group diverges from the classical geodesic-flow at a multiple of the Ehrenfest time--- namely, the former escapes to the cusp at logarithmic times, while the latter remains equidistributed on the surface.  Naturally, to show such a phenomenon we will have to analyze the wave propagation to large logarithmic times--- precisely what the current microlocal technology cannot do--- and so we turn to the number-theoretic tools available in the special case of Eisenstein series on a congruence surface.

\vspace{.2in}

We recall that the Eisenstein series are given by
$$E(z,s) = \sum_{\gamma\in\Gamma_\infty\backslash\Gamma} Im(\gamma.z)^{s}$$
for $Re(s)>1$, and analytically continued to be meromorphic in the entire plane.  We will be interested in the half-line $E(z,\frac{1}{2}+ir)$ as $r\to\infty$.
These are eigenfunctions are the Laplacian, of eigenvalue $\frac{1}{4} + r^2$, but are not in $L^2(M)$.  However, Luo-Sarnak \cite{LuoSarnak} showed that one can exploit the arithmetic structure of the surface to analyze expectation values for Eisenstein series explicitly.  In particular, they show that for fixed smooth functions $f, g \in C_0^\infty(M)$ we have
\begin{equation}\label{LS-Jak}	
\frac{\mu_r(f)}{\mu_r(g)} = \frac{\int_{M} f(x) dx}{\int_{M} g(x) dx}	
\end{equation}
where 
$$\mu_r(g) = \left\langle g\cdot E\left(z, \frac{1}{2}+ir\right), E\left(z,\frac{1}{2}+ir\right)\right\rangle$$
gives the expectation values for the state $E(z,\frac12+ir)$ of Eisenstein series.
The renormalization is important, as they actually show that\footnote{The original paper \cite{LuoSarnak} has some typos, leading to an incorrect constant; we thank Yiannis Petridis for pointing  out this correction, as it appears in the paper \cite{PetridisRR, PetridisRR_Err}.}
$$\mu_r(f) \sim \frac{6}{\pi}\log{r} \cdot \int_{S^*M} f(x) dx$$
grows proportionally to $\log{r}$.  Jakobson \cite{Jak} extended this result to microlocal lifts of $\mu_r$ to $S^*M$ as well.

It is important to note that such an explicit analysis of expectation values is only possible here because of the number-theoretic tools available in this context.  In particular, estimates for the expectation values of test functions $g$ orthogonal to constants are derived from subconvexity bounds for the Riemann zeta function and $L$-functions associated to Maass cusp forms.  On one hand, these powerful tools will allow us to analyze wave propagation of Eisenstein quasimodes to large multiples of Eherenfest time, but on the other hand we are restricted to this highly non-generic setting of congruence surfaces, with their special arithmetic symmetries.  

As remarked in \cite{EisenQuasi}, there are some very good reasons to believe that this behavior of Eisenstein quasimodes on congruence surfaces should not hold generically.  The rapid escape of these quasimodes to the cusp can be seen as akin to the concentration of quasimodes on geodesics observed in \cite{loc_quasi} for a compact surface, and since non-periodic geodesics on a compact surface do not escape but rather equidistribute, we could expect to see agreement between classical and quantum evolution persist in the compact setting.  But even in the non-compact case, we expect that this rapid escape is a special feature of the arithmetic structure, and tied to  the dearth of continuous spectrum on congruence surfaces as first observed by Selberg (see eg. \cite[11.1]{Iwaniec}).  This is expected to be a very special feature of the arithmetic structure, and generically the discrete spectrum should be small and the continuous spectrum dominant \cite{PhillipsSarnak}.  We expect therefore that generically, quasimodes in the continuous spectrum should not escape to the cusp until the Heisenberg time, in contrast to the phenomenon exhibited here.  Since we have virtually no tools to analyze spectral data on generic surfaces, it would be difficult to prove any statements in this direction; on the other hand, it might be interesting to look at numerics for, say, non-arithmetic Hecke triangle groups, and see if the semiclassical correspondence (\ref{sc c}) persists here to the longer times observed by \cite{TomsovicHeller,TomsovicHeller2} for billiards.

\vspace{.2in}
{\bf Acknowledgements.}  We thank Roman Schubert for many helpful discussions, and in particular for introducing us to the subject of long-time propagation of unstable Lagrangian states.

\section{Eisenstein Quasimodes}

We recall the following result of \cite{EisenQuasi} regarding ``Eisenstein quasimodes" on the modular surface $M=SL(2,\mathbb{Z})\backslash\mathbb{H}$.  Let $\{h_j\}$ be a sequence of smooth functions such that $\int_r h_j(r)dr = 1$, and satisfying
$$\int_r h_j(r) |r-r_j| = o(\log\log{r_j}/\log{r_j})$$
for a corresponding sequence of ``approximate spectral parameters" $r_j\to\infty$.  Define the Eisenstein quasimodes\footnote{We note that our definition here differs from a sign in the Eisenstein series; this is for convenience when we consider incoming plane waves of negative spectral parameter.  This changes some signs in the formulas of Theorem~\ref{sharp asymptotic} and Corollary~\ref{Ehrenfest cutoff}.}
$$E_{h_j}(z) = \int_r h_j(r) E\left(z, \frac{1}{2}-ir\right) dr$$
and the corresponding microlocal lift distributions
$$\mu_j(f) := \left\langle Op(f) E_{h_j}, E_{h_j}\right\rangle$$
In \cite{EisenQuasi} we studied weak-* limits of these distributions, and in particular derived the following asymptotic:
\begin{Theorem}\label{sharp asymptotic}
Let $\{h_j\}$ as above, and the microlocal lifts $\mu_j$ of the Eisenstein quasimodes $E_{h_j}$.  Then for any smooth, compactly supported test function $f\in C_0^\infty(S^*M)$, we have
$$\frac{\mu_{h_j}(f)}{2\log{r_j}} + o(1) = \frac{3}{\pi}\int_{S^*M} f(x)dx \cdot \iint_{r_1,r_2} h_j(r_1) \overline{h_j}(r_2) \frac{e^{2i\log{r_j}\cdot(r_2-r_1)}-1}{2i\log{r_j}\cdot (r_2-r_1)}dr_1dr_2 $$
\end{Theorem}

A simple consequence of Theorem~\ref{sharp asymptotic} is the following
\begin{Corollary}\label{Ehrenfest cutoff}
Let $\{h_j\}$ be a sequence of smooth, real-valued functions such that $\int_r h_j(r)dr=1$ and $\int_r h(r) |r-r_j|dr \lesssim 1/\log{r_j}$.  Then
 the corresponding microlocal lift $\mu_j$ is asymptotic to
$$\mu_j(f) \sim \frac{3}{\pi} \int_{S^*M}f(x)dx\cdot   \int_{t=0}^{2\log{r_j}}|\widehat{h_j}(t)|^2 dt$$
for every $f\in C_0^\infty(S^*M)$.
\end{Corollary}
Since the volume of the modular surface is $\frac{\pi}{3}$, the quantity $\frac{3}{\pi}\int f(x)dx$ is the normalized average of $f$ on $S^*M$; thus, in other words, Corollary~\ref{Ehrenfest cutoff} says that the total mass of $\mu_{h_j}$ that does not escape to the cusp is given by $\int_{t=0}^{2\log{r_j}}|\widehat{h_j}(t)|^2 dt$.

Since $\widehat{h_j}$ represents wave propagation times, we interpret this result as saying that the Eisenstein plane waves follow geodesic paths and equidistribute on the surface up to a logarithmic time $2\log{r_j}$, and then ``immediately afterwards" disengage from their geodesic paths and retreat back to the cusp.  This can be understood as a breakdown of the semiclassical correspondence of wave propagation to geodesic flow, at a fixed multiple of the Ehrenfest time; we recall that this semiclassical correspondence is known to hold up to the Ehrenfest time $T_E$, and one suspects that generically it may hold for much longer times, perhaps even to a power of $r$.

{\em Proof of Corollary~\ref{Ehrenfest cutoff}:  }
To prove the Corollary simply apply Plancherel and the usual properties of convolution and Fourier transform to get
\begin{eqnarray*}
\frac{\mu_{h_j}(f)}{2\log{r_j}} & \sim & \left(	\frac{3}{\pi}\int_{S^*M}f(x)dx	\right) \int_{r_1}{h_j}(r_1) {\int_{r_2} \overline{h_j}(r_2) \frac{e^{2i\log{r_j}\cdot (r_2-r_1)}-1}{2i\log{r_j}\cdot (r_2-r_1)}dr_2	} dr_1\\
& = &  \left(	\frac{3}{\pi}\int_{S^*M}f(x)dx	\right) \Big\langle {h_j}, {{h_j}\ast u(2\log{r_j}(\cdot)) 	}\Big\rangle\\
& = & \left(	\frac{3}{\pi}\int_{S^*M}f(x)dx	\right) \left\langle \widehat{{h_j}}, {\widehat{{h_j}}\cdot \frac{1}{2\log{r_j}}\chi([0, 2\log{r_j}])}\right\rangle
\end{eqnarray*}
since the Fourier transform of the function $u(s) = \frac{e^{is}-1}{is}$ is the characteristic function of the interval   $[0,1]$.  $\Box$
\vspace{.2in}

\section{Example of Semiclassical Breakdown}

In this section, we show how the above can be used to construct an example of a Lagrangian state associated to an unstable manifold, whose classical and quantum expectations diverge from each other at a multiple of the Ehrenfest time.  Our setting is the modular surface $M=SL(2,\mathbb{Z})\backslash \mathbb{H}$, and we will take our Lagrangian state to be invariant under the horizontal ($x$) coordinate, and compactly supported in the vertical ($y$) coordinate, inside the usual fundamental domain.  

Let $a\in C^\infty(\mathbb{R}^+)$ be fixed, positive, and supported in the interval $[2,3]$.  We consider for $-\frac12<x\leq \frac12$ the function 
$$f_{r}(x+iy) = a(y) y^{-ir}$$
projected to $M$ by summing over translations on the left by $\Gamma=SL(2,\mathbb{Z})$ (note the original function is compactly supported inside the fundamental domain),
which is a Lagrangian state associated to the horocycles $y=const.$ relative to the cusp at $\infty$, viewed as unstable horocycles since the direction given by the gradient of the phase function $-\log(y)$ is negative, pointing downward away from the cusp.

We now wish to study the propagation of these Lagrangian states for time $C\log{r}$, where $C$ is some large constant, and compare with the classical flow  along geodesics.  Since the states $f_{r}$ are independent of the $x$ variable and thus constant along the periodic horocycles relative to the cusp, they are orthogonal to Maass cusp forms,  and thus essentially supported on the continuous spectrum spanned by Eisenstein series; the constant function component is dealt with separately and 
will be negligible\footnote{We are using the fact that the modular surface has no exceptional eigenvalues, but the argument generalizes to other congruence surfaces, and  any hypothetical exceptional eigenfunctions are dealt with the same way as the constants.}.
We compute the spectral decompostion \cite[Theorem 7.3]{Iwaniec}
$$f_{r} = \hat{f}_{r}(0) + \frac{1}{4\pi} \int_{-\infty}^\infty \left\langle f_{r} , E(\frac{1}{2}+is)\right\rangle E(\frac{1}{2}+is)ds$$
where $\hat{f}_{r}(0) = \int_0^\infty a(y)y^{-ir-2}dy$ is the constant component, and $E(\frac{1}{2}+is)$ is the usual Eisenstein series
$$E(z,z') = \sum_{\gamma\in \Gamma_\infty\backslash\Gamma} Im(\gamma z)^{z'}$$
for $Re(z')>1$, analytically continued to the line $Re(z')=\frac12$.
Here, $\Gamma=SL(2,\mathbb{Z})$ and $\Gamma_\infty$ is the parabolic subgroup 
$$\Gamma_\infty = \left\{\begin{pmatrix} 1& n\\ & 1\end{pmatrix} : n\in\mathbb{Z} \right\}$$
stabilizing the cusp at $\infty$.  Recall that the Fourier expansion of $E(\frac{1}{2}+is)$ in the $x$-coordinate (relative to the cusp) has zero-th coefficient equal to $y^{\frac{1}{2}+is} + \phi(\frac{1}{2}+is)y^{\frac{1}{2}-is}$, and  all non-zero coefficients are again orthogonal to $f_{r}$ since the latter is constant in $x$ on the fundamental domain.  We recall that the scattering coefficient 
$$\phi(\frac{1}{2}+is) = \frac{\xi(is)}{\xi(is+1)}=\sqrt{\pi}\frac{\Gamma(i\frac{s}{2})\zeta(is)}{\Gamma(\frac{1}{2}+i\frac{s}{2})\zeta(1+is)}$$
has absolute value $1$, by the functional equation $\xi(s)=\xi(1-s)$ and the identity $\xi(\bar{s})=\overline{\xi(s)}$.

Thus since $a(y)$ is compactly supported inside the fundamental domain, we have
\begin{eqnarray*}
f_{r} & = & \hat{f}_{r}(0) + \frac{1}{4\pi} \int_{-\infty}^\infty \left(\int_M a(y) y^{-ir} [y^{\frac{1}{2}-is}	+\phi(\frac{1}{2}-is)y^{\frac{1}{2}+is}] \frac{dxdy}{y^2} \right) E(\frac{1}{2}+is) ds\\
& = & \hat{f}_{r}(0) + \frac{1}{4\pi} \int_{-\infty}^\infty \left(\int_0^\infty a(y)  [y^{-\frac{3}{2}-i(s+r)}	+\phi(\frac{1}{2}-is)y^{-\frac{3}{2}+i(s-r)}] dy \right) E(\frac{1}{2}+is) ds
\end{eqnarray*}
Writing $\hat{a}(\eta) = \int_{0}^\infty a(y) y^{-\frac{1}{2}+i\eta}\frac{dy}{y}$, for the Mellin transform of $a$ at $(-\frac{1}{2}+i\eta)$, we get
$$f_{r} = \hat{f}_{r}(0) + \frac{1}{4\pi} \int_{-\infty}^\infty [\hat{a}(-s-r) + \phi(\frac{1}{2}-is)\hat{a}(s-r)] E(\frac{1}{2}+is) ds$$
Since $a$ is fixed, smooth, and compactly supported, it is clear that $\hat{a}(\eta)=O(|\eta|^{-\infty})$ decays rapidly, i.e. that for every $N$ there exists a $C_N$ such that $\hat{a}(\eta) \leq C_N(|\eta|^{-N})$.  We also note similarly that $\hat{f}_{r}(0) = O(r^{-\infty})$.

We now consider the family of operators $V_{T}$ given by
\begin{equation}\label{V_T}
V_{T} f = \frac{1}{T} \int_0^T e^{it(\sqrt{-\Delta-1/4} + r)} f dt
\end{equation}
These operators average over the wave propagation of $f$ to time $ T$, which for our purposes will be taken to be a large multiple of Ehrenfest time, i.e. $T=C\log{r}$ for some sufficiently large constant $C$.  We adjust each time propagation by a scalar phase factor of $e^{itr}$ in order to localize $V_{T}f$ spectrally around spectral parameter $-r$.  Applying this to the Lagrangian state $f_{r}$ defined above, we first analyze the resulting function spectrally to determine that for $T=C\log{r}$
\begin{equation}\label{spectral lower bound}
||V_{T} f_{r}||_{L^2(M)}^2 \gtrsim \frac{1}{C\log{r}}
\end{equation}
noting first that $||V_{T}\hat{f}_{r}(0)|| = O(r^{-\infty})$, and then computing
\begin{eqnarray*}
\lefteqn{V_{T} f_{r}}\\
 & = & O(r^{-\infty})
  + \frac{1}{4\pi} \int_{-\infty}^\infty \left(\frac{1}{T} \int_0^{T}	[\hat{a}(-s-r) + \phi(\frac{1}{2}-is)\hat{a}(s-r)] E(\frac{1}{2}+is) \cdot e^{it(s+r)} dt 	\right) ds\\
& = & O(r^{-\infty}) + \frac{1}{4\pi}\int_{-\infty}^\infty \left(\frac{1}{T}\int_0^T e^{it(s+r)} dt\right)\hat{a}(-s-r) E(\frac{1}{2}+is)ds	\\
& = & O(r^{-\infty}) + \frac{1}{4\pi} \int_{-\infty}^\infty \frac{e^{iT(s+r)}-1}{iT(s+r)} \hat{a}(-s-r) E(\frac{1}{2}+is) ds 
\end{eqnarray*}
Here we used the fact that $\hat{a}(s-r)= O(r^{-\infty})$ for $s$ near $-r$, to deduce that $\frac{1}{T}\int_0^T e^{it(s+r)} \hat{a}(s-r)dt=O(r^{-\infty})$ as well, and incorporate it into the existing error term.

By the unitarity of the Eisenstein transform \cite[Proposiiton 7.1]{Iwaniec}, we see that since the real part of the coefficient $\frac{e^{iT(s+r)}-1}{iT(s+r)}$ is $\frac{\sin{T(s+r)}}{T(s+r)}$, we will get the estimate
\begin{eqnarray*}
||V_{T}f_{r}||_{L^2(M)}^2 & \gtrsim &   \int_{-\infty}^\infty  \frac{\sin^2(T(s+r))}{T^2(s+r)^2}|\hat{a}(-s-r)|^2 dr\\
& \gtrsim & \int_{|s+r|\leq T^{-1}} \frac{\sin^2(T(s+r))}{T^2(s+r)^2}|\hat{a}(-s-r)|^2 dr
\end{eqnarray*}
Since $a$ is positive, we deduce that $|\hat{a}(\eta)|$ is bounded away from $0$ on a neighborhood of $\eta=0$, and in particular $|\hat{a}(-s-r)|^2\gtrsim 1$ for $|s+r|\leq T^{-1}$.   Thus
since the factor $\frac{\sin^2(T(s+r))}{T^2(s+r)^2} \gtrsim 1$ for $|s+r|\leq T^{-1}$ as well, we deduce the desired bound (\ref{spectral lower bound}).

We note also that the above calculation shows that $V_{T}f_{r}$ is an Eisenstein $O(1/\log{r})$ quasimode (recall we set $T=C\log{r}$), satisfying the hypotheses of Corollary~\ref{Ehrenfest cutoff}.  The Fourier transform of the spectral coefficients $h_r(s) = \frac{e^{iT(s+r)}}{iT(s+r)} \hat{a}(-s-r)$ of $V_{T}f_{r}$ is everywhere bounded by $O(T^{-1})$--- eg., since the Fourier transform of $\frac{e^{iT(s+r)}}{iT(s+r)}$ is everywhere bounded by $T^{-1}$, and $|\hat{a}|$ is fixed and its integral doesn't depend on $T$---  and so for any compactly supported test function $g\in C_0^\infty(M)$ we have by Corollary~\ref{Ehrenfest cutoff}
$$\mu_{h_j}(g) \lesssim \int_{t=0}^{ 2\log{r}} |O(T^{-1})|^2dt \lesssim \log{r} \cdot (C\log{r})^{-2} = O(1/C^2\log{r})$$
Since by (\ref{spectral lower bound}) the full $L^2$-norm of $V_{T}f_{r}$ is $\gtrsim 1/C\log{r}$, we see that for $C$ sufficiently large, an arbitrarily high percentage of the $L^2$-mass of $V_{T}f_{r}$ is concentrated in the cusp.  For we may take $g$ to be identically $1$ on an arbitrarily large compact subset of $M$, so that the above estimate says we may concentrate an arbitrarily large percentage of the mass in an arbitrarily small neighborhood of the cusp, by taking $C$ large enough.  It is from this property that we will derive a breakdown of the semiclassical correspondence for large multiples of Ehrenfest time.

To complete the argument, we consider the classical expectations of our initial state, pushed-forward by the geodesic flow.  Since our initial Lagrangian state is $x$-independent and supported in $2\leq y\leq 3$, this can be estimated using the quantitative equidistribution of long horocycles due to Sarnak:
\begin{Theorem}[\cite{SarnakHoro}]\label{horo}
Consider a long periodic horocycle $\eta_y = \{x+iy : -\frac{1}{2}\leq x\leq \frac{1}{2} \}$ as $y\to 0$; note that $\eta$ has length $y^{-1}$.  Then there exists $\delta>0$ such that for any compactly supported test function $g\in C_0^\infty(M)$, we have
$$y \int_\eta g(x+iy) \frac{dx}{y} = \int_M g dArea + O(y^\delta)$$
\end{Theorem}

{\em Remark:}  In fact, Sarnak \cite{SarnakHoro} showed much more, that one can replace the error with a sharper $o(y^{1/2})$, but we will not need the full strength of this result.

{\em Completion of the proof:}  
We have shown that the wave-propagation-averaged functions $V_Tf_r$ localize in the cusp, in the sense that
$$||V_Tf_r||_2^2 \gtrsim 1/C\log{r}$$
while for any function $g$ of compact support we have
$$\langle g V_Tf_r, V_Tf_r\rangle \lesssim 1/C^2\log{r}$$
We now wish to contrast this with the push-forward of $f_r$ by the classical geodesic flow, i.e. the classical expectation
$$\langle g_{\Lambda(T)} W_Tf_r, W_Tf_r\rangle$$
where $g_\Lambda(T) = g\cdot \delta_{\Lambda(T)}$ is the restriction of $g$ to the  Lagrangian union of horocycles $\Lambda(T) = \bigcup_{y=2e^{-T}}^3 \eta_{y}$.  The operator $W_T$ is the classical analog of the quantum operator $V_T$, and is given here by averaging the unitary geodesic flow along $\Lambda(T)$ (multiplied by the appropriate phases as in (\ref{V_T}) for $V_T$)
\begin{eqnarray*}
W_Tf_r (y) & = & \frac{1}{T} \int_{t=0}^T f_r(ye^{t})e^{-t/2}e^{itr}dt\\
& = & \frac{1}{T} \int_{t=0}^T a(ye^{t})\cdot [ye^{t}]^{-ir} e^{itr} e^{-t/2} dt\\
& = & \frac{1}{T} \int_{t=0}^T a(ye^{t})e^{-t/2} y^{-ir}dt 	
\end{eqnarray*}
For simplicity, we take $g\in C^\infty(M)$ to depend only position, constant on the $SO(2)$-fibres in $S^*M$.

Since the lift of the Lagrangian to $S^*\mathbb{H}$ is given by unstable horocycles relative to the cusp, and the geodesic flow maps this family of horocycles into itself, we can easily show that this mapping is one-to-one on $S^*M$ for arbitrary times $T$.  For if a closed geodesic mapped a point on one of these horocycles to a point on another, it would require a non-trivial element  $\gamma\in \Gamma_\infty\backslash\Gamma$ to lie in the Borel subgroup $\gamma\in AN = \left\{\begin{pmatrix} a & b\\ 0 & d\end{pmatrix}\right\}$, but the requirement that $\gamma\in SL(2,\mathbb{Z})$ would then mean that $a=d=\pm 1$, and the only such elements in $SL(2,\mathbb{Z})$ belong to  the parabolic subgroup $\Gamma_\infty$.  Thus we can 
compute the classical expectation as
\begin{eqnarray*}
\langle g_{\Lambda(T)} W_Tf_r, W_Tf_r\rangle & = & \int_{y=0}^\infty \int_{x=0}^1 g(x+iy) |W_Tf_r(y)|^2 \frac{dxdy}{y^2}\\
 & = & \int_{y=0}^\infty |W_Tf_r(y)|^2 \int_{x=0}^1 g(x+iy) \frac{dxdy}{y^2}\\
& = & \frac{1}{T^2}\int_{y=0}^\infty \left( \int_{t=0}^T a(ye^{t}) e^{-t/2} dt\right)^2 \int_{x=0}^1 g(x+iy) \frac{dxdy}{y^2}\\
& = & \frac{1}{T^2}\int_{y=0}^\infty \left( \int_{u=y}^{ye^T} a(u) \sqrt{\frac{y}{u}} \frac{du}{u}\right)^2 \int_{x=0}^1 g(x+iy) \frac{dxdy}{y^2}\\
& = & \frac{1}{T^2}C_a^2 \int_{y=2e^{-T}}^3 \left[y \int_{\eta_y} g(x+iy) \frac{dx}{y}\right] \frac{dy}{y} + O(T^{-2})
\end{eqnarray*}
where $C_a = \int_0^\infty a(u)u^{-3/2}du$ is a bounded constant,
since $a$ is fixed, positive and compactly supported.  

 Using the quantitative equidistribution of horocycles (Theorem~\ref{horo}), 
\begin{eqnarray*}
\langle g_{\Lambda(T)} W_Tf_r, W_Tf_r\rangle & = & \frac{1}{T^2} C_a^2 \int_{2e^{-T}}^3 \left( \int_M g \cdot dArea + y^\delta\right) \frac{dy}{y} +O(T^{-2})\\
& = & \frac{1}{T}C_a^2 \int_M g \cdot dArea + O(T^{-2})\\
& \gtrsim & \frac{1}{T} \gtrsim \frac{1}{C\log{r}}
\end{eqnarray*} 
when $T=C\log{r}$.

So combining with (\ref{spectral lower bound}),  we have finally
\begin{eqnarray*}
\langle g V_Tf_r, V_Tf_r\rangle & \lesssim & 1/C^2\log{r}\\
\langle g_{\Lambda(T)} W_Tf_r, W_Tf_r\rangle & \gtrsim &  1/C\log{r}
\end{eqnarray*}
so that for $C$ sufficiently large--- i.e., for $T$ a large multiple of Ehrenfest time--- the quantum expectation computed from $V_Tf_r$ diverges from the classical expectation given via $W_Tf_r$.
$\Box$

\def\cprime{$'$}

\end{document}